\documentclass{paper}
\usepackage[dvips]{color}
\usepackage{graphics,graphicx}
\usepackage{latexsym}
\usepackage{times}
\usepackage{natbib} 
\usepackage{amssymb,amsmath,mathrsfs,amsfonts,amsthm}
\usepackage{amsbsy} 
\usepackage{enumerate}
\newtheorem{theorem}{Theorem}[section]

\newtheorem{lemma}{Lemma}[section]

\usepackage{color}

\def\bR{\mathbb{R}}

\newcommand{\R}{\mathbb{R}} 

\newcommand{\bmat}{\begin{pmatrix}}
\newcommand{\emat}{\end{pmatrix}}
\newcommand{\beq}{\begin{equation}}
\newcommand{\eeq}{\end{equation}}

\usepackage{color}

\newcommand{\diag}{\mathrm{diag}}

\newcommand{\tr}{\mathrm{tr}}
\newcommand{\tran}{\mathrm{T}}

\newcommand{\mA}{\mathcal{A}}

\begin{document}

\title{
\emph{\Large Shrinking the Sample Covariance Matrix using Convex Penalties on the Matrix-Log
Transformation} 
}

\author{\large David E. Tyler*  ~ {\large and} ~ Mengxi Yi  
\thanks{Research for both authors was supported in part by NSF Grant DMS-1407751.}\\
\normalsize{Department of Statistics \& Biostatistics} \\[-4pt]
\normalsize{Rutgers, The State University of New Jersey, Piscataway, NJ 08854, U.S.A.} \\[-4pt]
\normalsize{dtyler@stat.rutgers.edu; mengxiyi@stat.rutgers.edu}} 

\maketitle
\begin{abstract}
For $q$-dimensional data, penalized versions of the sample covariance matrix are important when the sample size is
small or modest relative to $q$. Since the negative log-likelihood under multivariate normal sampling
is convex in $\Sigma^{-1}$, the inverse of its covariance matrix, it is common to add to it a penalty which is also
convex in $\Sigma^{-1}$. More recently, \cite{Deng-Tsui:2013} and \cite{Yu:2017} have proposed penalties which are
functions of the eigenvalues of $\Sigma$, and are convex in $\log \Sigma$, but not in $\Sigma^{-1}$. The resulting
penalized optimization problem is not convex in either $\log \Sigma$ or $\Sigma^{-1}$. In this paper, we note 
that this optimization problem is geodesically convex in $\Sigma$, which allows us to establish the existence 
and uniqueness of the corresponding penalized covariance matrices. More generally, we show the equivalence of 
convexity in $\log \Sigma$ and geodesic convexity for penalties on $\Sigma$ which are strictly functions of their 
eigenvalues. In addition, when using such penalties, we show that the resulting optimization problem reduces to
to a $q$-dimensional convex optimization problem on the eigenvalues of $\Sigma$, which can then be readily
solved via Newton-Raphson. Finally, we argue that it is better to apply these penalties to the shape matrix
$\Sigma/(\det \Sigma)^{1/q}$ rather than to $\Sigma$ itself. A simulation study and an example illustrate the
advantages of applying the penalty to the shape matrix.
\end{abstract}

\noindent
\textbf{Keywords:} geodesic convexity; M-estimation; Newton-Raphson algorithm; penalized covariance matrices.

\section{Introduction and Motivation} \label{Sec:Intro}
For a $q$ dimensional sample $x_1, \ldots, x_n$, the sample covariance matrix 
$S_n = n^{-1}\sum_{i=1}^n (x_i - \bar{x})(x_i - \bar{x})^\tran$ is not well-conditioned and can be highly variable when $q$ is of the same order as $n$. In such cases, 
one may wish to consider a regularized or a penalized version of the sample covariance matrix.  Since the loss function obtained from the negative log likelihood under multivariate normal sampling 
\begin{equation} \label{eq:lik}
 l(\Sigma;S_n) = \tr(\Sigma^{-1}S_n) + \log \det \Sigma,
\end{equation} 
is convex in $\Sigma^{-1}$, it is natural to consider additive penalties which are also convex in $\Sigma^{-1}$ such as the 
the graphical lasso penalty $\sum_{i \ne j} | \{\Sigma^{-1}\}_{ij}|$ \ \citep{Yuan-Lin:2007,Friedman:2008}.  Minimizing the penalized loss function
\begin{equation} \label{eq:pen}
 L(\Sigma;S_n,\eta) = l(\Sigma;S_n) + \eta \Pi(\Sigma),
\end{equation}
over the set of symmetric positive definite matrices $\Sigma > 0$, with $\Pi(\Sigma)$ being a non-negative penalty function and $\eta \ge 0$ being a tuning parameter, is then a convex optimization problem. 

More recently, \cite{Deng-Tsui:2013} consider the penalty $\Pi_R(\Sigma) \equiv \|\log \Sigma \|_F$, where the norm 
refers to the Frobenius norm. This penalty is strictly convex in $\log \Sigma$ but not in $\Sigma^{-1}$. By letting
$A = \log \Sigma$, they observe that, when using this penalty, \eqref{eq:pen} can be express in terms of a penalized loss function 
over the set of symmetric matrices of order $q$, namely 
\begin{equation} \label{eq:Rem}
\mathcal{L}(A;S_n,\eta) = \tr(e^{-A}S_n) + \tr A + \eta \| A \|^2_F,
\end{equation}
with the penalty $\| A \|^2_F = \tr(A^2)$ being strictly convex in $A$. As noted in section \ref{Sec:G-L}, the function $\tr(e^{-A}S_n)$ is not 
in general a convex function of $A$, and consequently minimizing \eqref{eq:Rem} over $A$ does not correspond to a convex optimization problem. 
Hence, there is no assurance as to the existence and uniqueness of a minimum to \eqref{eq:Rem}. 

One of our objectives in this paper is to argue that rather than using the concept of convexity in $\log \Sigma$ in
problem \eqref{eq:Rem}, a more appropriate setting is based on the notion of geodesic convexity, or g-convexity for short.   
The function $\|\log \Sigma \|_F$  has been well studied within Riemannian geometry and corresponds to the Riemannian or geodesic distance 
between $\Sigma$ and the identity matrix \citep{Moakher:2005,Bhatia:2009}, and is known to be
strictly g-convex in $\Sigma$. For $S_n \ne 0$, 
the loss function \eqref{eq:lik} is also strictly g-convex, and consequently the penalized loss function \eqref{eq:pen} 
is strictly g-convex, when choosing $\Pi = \Pi_R$.  Moreover, \eqref{eq:pen}, with $\Pi = \Pi_R$, can be 
shown to be g-coercive, which implies it has a unique critical point, with this unique critical point 
corresponding to its global minimum; see Lemmas \ref{Lem:nonsing} and \ref{Lem:sing}.

The concept of g-convexity can be mathematically challenging, and in practice it can be difficult to prove that a given function is g-convex.  A further contribution of this paper is to show that for an orthogonally invariant penalties, i.e.\ penalties which are strictly functions of the eigenvalues of $\Sigma$, as is the case for $\Pi_R$, (strict) g-convexity in $\Sigma$ and (strict) convexity in $\log \Sigma$ are equivalent. Furthermore, it is shown that g-convexity for such function reduces to the simpler task of establishing 
(strict) convexity when viewed as function on the logs of the eigenvalues; see Theorem \ref{Thm:g-log}.  
For example, if we express $\Pi_R(\Sigma) = \sum_{j=1}^q a_j^2 $, where $a_j = \log \lambda_j$ with $\lambda_1 \ge \cdots \ge \lambda_q > 0$ being the eigenvalues of $\Sigma$, then it is strictly convex as a function of $a \in \R^q$, and hence $\Pi_R(\Sigma)$ is strictly convex in
$\log \Sigma$ and strictly g-convex in $\Sigma$.

\cite{Deng-Tsui:2013} also propose an iterative quadratic programming algorithm over the class of symmetric matrices $A$ 
of order $q$ for finding the minimum of \eqref{eq:Rem}.  We show in Theorem \ref{Thm:diag}, though, 
that the solution to this problem has the same eigenvectors as $S_n$. This leads to a simpler algorithm 
based on finding the minimum of a strictly convex univariate function for each eigenvalue, namely 
$g(a;d) = d e^{-a} + a - \eta a^2$, with $d$ corresponding to an eigenvalue of $S_n$ and $a$ being the corresponding log 
eigenvalue of $\Sigma$. The solution to this univariate convex optimization problem can be readily obtained via a 
Newton-Raphson algorithm.

As recently noted by \cite{Yu:2017}, the penalty $\| A \|^2_F$ shrinks the sample covariance matrix towards the 
identity matrix. They proposed using the alternative penalty $\| A - \widehat{m} I_q \|^2_F$, with $\widehat{m}$ being an 
estimate of the mean of the log of the eigenvalues of $\Sigma$, i.e.\ of  $m(A)  = \tr{A}/q$. Since $\widehat{m}$ is first 
determined from the data, this does not correspond to a pure penalty function for $\Sigma$. Rather than using a preliminary 
estimate of $m$, we propose replacing $\widehat{m}$ with $m(A)$. This approach yields an estimate of $m(A)$ consistent with the 
penalized estimate of $\Sigma$, i.e.\ $m(\widehat{A}) = \tr{\widehat{A}}/q$. The resulting penalized objective function 
\eqref{eq:pen}, when using the penalty $\| A - m(A) I_q \|^2_F = \sum_{j=1}^q (a_j-\overline{a})^2$ is shown, within 
section \ref{Sec:shape}, to also be strictly g-convex and g-coercive. Consequently, the global minimum of \eqref{eq:pen} 
corresponds to the unique critical point. The solution to this optimization problem reduces to finding the minimum of 
a strictly convex function in $\R^q$ .

Summarizing, this article is organized as follows. In section \ref{Sec:G-L}, the concept of geodesic convexity is briefly
reviewed, and results on the existence and uniqueness of penalized sample covariance
matrices based on g-convex penalty functions in general are presented. Results on the relationship between convexity 
in $\log \Sigma$ and g-convexity in $\Sigma$ are given in section \ref{Sec:Log}. In section \ref{Sec:shape}, convexity
results for \eqref{eq:pen} are given when applying the penalty to the shape matrix $\Sigma/(\det \Sigma)^{1/q}$ rather 
than to $\Sigma$ itself, with $\| A - m(A) I_q \|^2_F$ being a special case of such a shape penalty. Algorithms for 
computing the penalized sample covariance matrices, based on orthogonally invariant g-convex penalties are given in
section \ref{Sec:Opt}. We emphasize that this paper treats g-convex penalties in general, with applications to 
$\Pi_R$ treated as a special case. The results of a simulation study discussed in section \ref{Sec:Sim}, together with
an example given in section \ref{Sec:Ex}, demonstrate the advantages of penalizing shape.  Proofs and some technical
details are given in an appendix..

\section{Geodesic Convexity} \label{Sec:G-L}
The notion of geodesic distance between multivariate normal distributions, or equivalently the geodesic distance between their covariance
matrices, has been a topic of interest at least as early as \cite{Skovgaard:1984}. However, the realization the multivariate normal negative 
log-likelihood $l(\Sigma;S_n)$ is g-convex, and strictly g-convex when $S_n \ne 0$, which follows as a 
special case of Theorem 1 in \cite{Zhang:2013}, is relatively recent.

The set of symmetric positive definite matrices of order $q$ can be viewed as a Riemannian manifold with the 
geodesic path from $\Sigma_0 > 0$ to $\Sigma_1 > 0$ being given by 
$\Sigma_t = \Sigma_0^{1/2} \{\Sigma_0^{-1/2}\Sigma_1\Sigma_0^{-1/2}\}^t \Sigma_0^{1/2}$ for $0 \le t \le 1$, 
see \cite{Bhatia:2009} or \cite{Wiesel-Zhang:2015} for more details. An alternative representation 
for this path is given by $\Sigma_t = Be^{t\Delta} B^\tran$, where $\Sigma_0 = BB^\tran$ and $\Sigma_1 = Be^\Delta B^\tran$ 
with $\Delta$ being a diagonal matrix of order $q$. A function $f(\Sigma)$ is said to be g-convex if and only if 
$f(\Sigma_t) \le (1-t)f(\Sigma_0) + tf(\Sigma_1)$ for $0 < t < 1$, and it is strictly g-convex if strict inequality holds 
for $\Sigma_0 \ne \Sigma_1$. Analogous to convexity in $\log \Sigma$, for which convexity in $\log \Sigma$ implies 
convexity in $\log \Sigma^{-1}$, g-convexity in $\Sigma$ implies g-convexity in $\Sigma^{-1}$.

As with convexity, any local minimum of a g-convex function is a global minimum, and when differentiable any critical 
point is a global minimum, with the set of all minima being g-convex. In addition, if a minimum exists, then the minimum 
is unique when the function is strictly g-convex. Finally, the sum of two g-convex functions is g-convex, and the sum is 
strictly g-convex if either of the two g-convex summands is strictly g-convex. Consequently, the following lemma holds.
\begin{lemma} \label{Lem:Lgeo}
If $\Pi(\Sigma)$ is g-convex and $S_n \ne 0$, then $L(\Sigma;S_n,\eta)$ is strictly g-convex on $\Sigma > 0$, and  
the set of all local minima $\mA_\eta$ is either empty or contains a single element. That is,
if there exists a minimizer $\widehat{\Sigma}_\eta > 0$ to $L(\Sigma;S_n,\eta)$, then it is unique.
\end{lemma}
The existence of a minima for a g-convex function requires some additional conditions, with a necessary and sufficient condition being that it be 
geodesic coercive \citep{Duembgen-Tyler:2016}. A g-convex function $F(\Sigma)$ is said to be g-coercive if and only if  
$F(\Sigma) \to \infty$ as $\| \log \Sigma \|_F \to \infty$. For $S_n > 0$, $l(\Sigma;S_n)$ is g-coercive and so, since
$\Pi(\Sigma)$ is bounded below, $L(\Sigma;S_n,\eta)$ is g-coercive and hence has a unique minimizer. Moreover, since a g-convex function
is continuous on $\Sigma > 0$, it follows that the solution is a continuous function of $\eta \ge 0$. This is summarized in the following lemma.
\begin{lemma} \label{Lem:nonsing}
Under the conditions of Lemma \ref{Lem:Lgeo}, if $S_n > 0$, then there exists a unique critical point $\widehat{\Sigma}_\eta > 0 $ to 
$L(\Sigma;S_n,\eta)$, with $\widehat{\Sigma}_\eta$ being its unique minimizer. Furthermore, $\widehat{\Sigma}_\eta$ is a continuous function of $\eta \ge 0$.
\end{lemma}

For singular $S_n$, some conditions on the penalty function are needed since it is possible for $\tr(\Sigma^{-1}S_n)$ to 
be bounded as $\log \det \Sigma  \to -\infty$, and hence $l(\Sigma;S_n)$ is not g-coercive in this case. A sufficient 
condition for $L(\Sigma,S_n,\eta)$ to be g-coercive when $S_n$ singular is that $\Pi(\Sigma)$ be g-coercive and $\eta > 0$.  
This condition, however, is too strong, and does not hold for the scale invariant or shape penalties discussed in section 
\ref{Sec:shape}. Some weaker conditions are given in the following lemma, with these conditions holding when $\Pi(\Sigma)$ is g-coercive. Note that under each of the three conditions
below, $\Pi(\Sigma) \to \infty$. 
\begin{lemma} \label{Lem:sing}
Under the conditions of Lemma \ref{Lem:Lgeo}, if  \\[2pt]
\begin{tabular}{rl}
(i)& $\Pi(\Sigma) \to \infty$ whenever $|\log \det \Sigma|$ is bounded above and $\| \log \Sigma \|_F \to \infty$, \\
(ii)& $(\log \det \Sigma)/\Pi(\Sigma) \to 0$ whenever $\log \det \Sigma \to -\infty$ but with $\lambda_1$ bounded away from $0$, and \\
(iii)& $\{\log(\lambda_1/\lambda_q)\}/\Pi(\Sigma)$ is bounded above whenever $\lambda_1 \to 0$ but with $\lambda_1/\lambda_q$ bounded away from $1$.
\end{tabular} \\[2pt]
then the conclusions stated in Lemma \ref{Lem:nonsing} hold when $S_n \ne 0$ is singular and $\eta > 0$.
\end{lemma}

\section{Geodesic Convexity and Convexity in Log} \label{Sec:Log}
In the following, we show that for orthogonally invariant functions, 
g-convexity in $\Sigma$ is equivalent to convexity in $\log \Sigma$. We say that a function $F$ on the set of positive definite
matrices is orthogonally invariant if and only if $F(\Sigma) = F(H\Sigma H^\tran)$ for any orthogonal matrix $H$ of order $q$.
It is straightforward to show that such functions can then be expressed in terms of a symmetric function of its eigenvalues 
$\lambda_1 \ge \cdots \ge \lambda_q > 0$.
\begin{lemma} \label{Lem:orth}
The function $F(\Sigma)$ is orthogonally invariant if and only if for some symmetric, i.e. permutation invariant, function 
$f: \bR^q\to\bR, F(\Sigma) \equiv f(a_1, \ldots, a_q)$ where $a_j = \log \lambda_j, \ j = 1, \ldots, q$.
\end{lemma}
\begin{theorem} \label{Thm:g-log}
For an orthogonally invariant function $F(\Sigma)$, the following three conditions are equivalent: \\[2pt]
\begin{tabular}{rl}
(i)& $F(\Sigma)$ is (strictly) g-convex. \\
(ii)& $F(\Sigma)$ is (strictly) convex in $\log(\Sigma)$. \\
(iii)& The corresponding function $f$, as defined in lemma \ref{Lem:orth}, is (strictly) convex.
\end{tabular}
\end{theorem}

A clarifying point regarding Theorem \ref{Thm:g-log} may be helpful. It should be noted, for example, that the corresponding function on $\R^q$
for the log concave function $F(\Sigma) = \log \lambda_q$ is not $f(a_1, \ldots, a_q) = a_q$ which is linear, hence convex, but not symmetric. 
Rather, its corresponding function is $f(a_1, \ldots, a_q) = \min\{a_1, \ldots, a_q\}$ which is symmetric but concave.

Outside of orthogonal invariant functions, g-convexity and convexity in log do not necessarily coincide. For example, as previously noted,
$l(\Sigma;S_n)$ is strictly g-convex, but not necessarily convex in log.  In particular, although $\log \det \Sigma = \tr(A)$ is linear and hence
convex in $A$, whether or not the term $\tr(\Sigma^{-1}S_n) = \tr(e^{-A}S_n)$ is convex in $A$ depends on the value of $S_n$. For example, when $S_n = I$, the convexity of $\tr(e^{-A})$ follows from Theorem \ref{Thm:g-log} since $\sum_{j=1}^q e^{-a_j}$ is convex. As far as we are
aware, general conditions on $S_n$ needed for $\tr(e^{-A}S_n)$ to be convex have not been formally addressed in the literature.  
An example of $S_n$ for which $\tr(e^{-A}S_n)$ is not convex in $A$ is given in the appendix. On the other hand, an example of a function 
which is convex in log but not g-convex is also presented in the appendix.

We now apply these results to the penalty studied by \cite{Deng-Tsui:2013}, i.e.\ $\Pi_R(\Sigma) = \| \log \Sigma \|_F^2$.
This penalty is orthogonally invariant and can be expressed as $\Pi_R(\Sigma) = \sum_{j=1}^q a_i^2$, which is symmetric and strictly convex as a 
function of $a \in \R^q$. Hence, by Theorem \ref{Thm:g-log}, $\Pi_R$ is strictly g-convex, and so Lemma \ref{Lem:nonsing} holds. 
Furthermore, Lemma \ref{Lem:sing} also holds since $\log \det \Sigma/\Pi_R(\Sigma) = \{\sum_{j=1}^q a_i\}/\{\sum_{j=1}^q a_i^2\} \to 0$ 
as $\sum_{j=1}^q a_i \to -\infty$.

The geodesic convexity of $\| \log \Sigma \|_F^2$ has been previously established using more involved proofs, see \cite{Bhatia:2009} 
for comparison. The importance of Theorem \ref{Thm:g-log} is that for penalty functions which are strictly functions of the eigenvalues of $\Sigma$, 
it completely characterizes g-convexity, as well as provides a simple condition for verifying g-convexity. For example, it 
readily follows that the Kullback-Leibler divergence from the identity matrix, i.e.\ $\tr(\Sigma^{-1}) + \log \det(\Sigma)$, 
which is convex in $\Sigma^{-1}$ is also g-convex. The condition number penalty $\lambda_1/\lambda_q$ and the penalty  
$\tr(\Sigma)+\tr(\Sigma^{-1})$, among others considered by \cite{Wiesel:2012} and \cite{Duembgen-Tyler:2016}, are also 
seen to be g-convex.

\section{Penalizing the shape matrix} \label{Sec:shape}

Any penalty on $\Sigma > 0$ can also be applied to its shape matrix $V(\Sigma) = \Sigma/\det(\Sigma)^{1/q}$. 
Here $\det V(\Sigma) = 1$, with the orbits of $V(\Sigma)$ form equivalence classes over $\Sigma > 0$ \cite{Davy:2008}. 
This then generates the new penalty $\Pi_s(\Sigma) \equiv \Pi(V(\Sigma))$. If the original penalty is minimized e.g.\  at $\Sigma = I$,
then the new penalty is minimized at any $\Sigma \propto I$.  Applying the penalty studied by \cite{Deng-Tsui:2013} to the
shape matrix yields
\[\Pi_{R,s}(\Sigma) \equiv \Pi_R\{V(\Sigma)\} = \| \log \Sigma - q^{-1}\{\log \det \Sigma\}I \|^2_F = \| A - m I \|^2_F,\] 
where $m = q^{-1}\{\log \det \Sigma\}  = \tr{A}/q$.
Since $\Pi_{R,s}$ is orthogonally invariant, with $\Pi_{R,s}(\Sigma) = $ \mbox{$\sum_{j=1}^q (a_i - \overline{a})^2$} being convex, it follows 
from Theorem \ref{Thm:g-log} that $\Pi_{R,s}$ is convex in log as well as g-convex, although the convexity is not strict in this case. 
Thus, for non-singular $S_n$, Lemma \ref{Lem:nonsing} on existence and uniqueness applies when using $\Pi_{R,s}$ as the penalty term.  
Also, as shown in the appendix, the additional conditions given in Lemma \ref{Lem:sing} needed to assure existence and uniqueness when
$S_n$ is singular also holds when using this penalty.

More generally, applying any g-convex penalty or penalty which is convex in $\log \Sigma$ to the shape matrix of $\Sigma$, 
yields respectively a new g-convex penalty or penalty convex in $\log \Sigma$. The following theorem applies to any such
penalties and does not presume $\Pi$ is orthogonally invariant.
\begin{theorem} \label{Thm:shape} \mbox{  } \\[2pt]
\begin{tabular}{rl} 
(i) & If $\Pi(\Sigma)$ is g-convex, then $\Pi_s(\Sigma)$ is also g-convex. \\
(ii) & If $\Pi(\Sigma)$ is convex in $\log \Sigma$, then $\Pi_s(\Sigma)$ is also convex in $\log \Sigma$. 
\end{tabular}
\end{theorem}

Thus, for g-convex $\Pi(\Sigma)$, Lemma \ref{Lem:nonsing} on existence and uniqueness for the case when $S_n$ is non-singular still applies 
when the penalty term $\Pi$ is replaced by $\Pi_{s}$. 
As another example, if we apply the Kullback-Leibler divergence from the identity to the shape matrix of $\Sigma$, 
one obtains the penalty $\Pi_s(\Sigma) = \tr\{V(\Sigma)^{-1}\} = q \overline{\lambda}_g/\overline{\lambda}_h$, where
$\overline{\lambda}_g$ and $\overline{\lambda}_h$ are respectively the geometric mean and the harmonic mean of the 
eigenvalues of $\Sigma$. This ratio represents a measure of eccentricity for $\Sigma$,
and is minimized at any $\Sigma \propto I$. By the previous theorem, this new penalty is also g-convex, and hence 
Lemma \ref{Lem:nonsing} applies. It can be verified that Lemma \ref{Lem:sing} also applies for this case.

\section{Optimizing the penalized loss function} \label{Sec:Opt}
As noted in the introduction, \cite{Deng-Tsui:2013} propose a quadratic iterative programming algorithm over 
$A = \log \Sigma$. The algorithm is derived by a repeated application of the  Volterra integral equation for 
$e^{tA}$ to obtain a second order expansion. Although they state in their introduction that some other previously 
proposed ``methods have retained the use of the eigenvectors of $S_n$ in estimating $\Sigma$ or $\Sigma^{-1}$,'' 
it is not clear if they recognize that the minimum  $\widehat{A}_\eta$ to \eqref{eq:Rem}, and hence 
$\widehat{\Sigma}_\eta = e^{\widehat{A}_\eta}$, also retain the same eigenvectors as $S_n$. As shown in the 
following theorem, this is true for any orthogonally invariant penalty.
\begin{theorem} \label{Thm:diag}  
Suppose $\Pi(\Sigma)$ is orthogonally invariant. Using the spectral value decomposition, express
$S_n = P_nD_nP_n^\tran$ with $P_n$ being an orthogonal matrix of order $q$, and where $D_n = \diag\{d_1, \ldots, d_q\}$. Then 
\[L(\Sigma; S_n,\eta) \ge L(P_n \Lambda P_n^\tran,S_n,\eta),\] where $\Lambda = \diag\{\lambda_1, \ldots, \lambda_q \}$. 
\end{theorem}
This lemma then implies that for orthogonally invariant penalties, the penalized covariance matrix has the form
$\widehat{\Sigma}_\eta = P_n D_\eta P_n^\tran$ for some diagonal matrix $D_\eta$. In particular, 
$D_\eta = \diag\{e^{\widehat{a}_{\eta,1}}, \ldots, e^{\widehat{a}_{\eta,q}}\}$ with $\widehat{a}_\eta$
being the minimizer over $a \in \R^q$ of
\begin{equation} \label{eq:fa}
L_q(a;d;\eta) = \sum_{j=1}^q d_je^{-a_j} + a_j + \eta \pi(a_1, \ldots, a_q).
\end{equation}
Here $\pi$ is the function on $\R^q$ corresponding to the function $\Pi$ as defined in Lemma \ref{Lem:orth}, and $d_1 \ge \cdots \ge d_q$ are
the eigenvalues of $S_n$. The function $L_q(a;d;\eta)$ is strictly convex whenever $\pi(a)$ is convex in $a \in \R^q$, which by Theorem \ref{Thm:g-log}
holds whenever $\Pi(\Sigma)$ is g-convex, or equivalently convex in $\log \Sigma$. Thus, for orthogonally invariant g-convex penalties, the minimization problem \eqref{eq:pen} reduces to the simpler and numerically well studied problem of minimizing a strictly convex function over $\R^q$.

The penalties proposed by \cite{Deng-Tsui:2013} and \cite{Yu:2017} are both of the form $\Pi(\Sigma) = \| A - c I_q \|_F^2$
with $c$ not dependent on $\Sigma$. For these cases, we have $\pi(a) = \sum_{j=1}^q (a_j - c)^2$ and so
$L_q(a;d;\eta) = \sum_{j=1}^q \{ d_je^{-a_j} + a_j + \eta (a_j - c)^2 \}$. Rather than using their proposed quadratic iterative 
programing algorithm over the set of symmetric matrices of order $q$ for this problem, one only needs to solve $q$ univariate 
strictly convex optimization problems, namely $\min\{ d_j e^{-a_j} + a_j + \eta (a_j -c)^2\}$ for $j = 1, \ldots, q$. 
The Newton-Raphson algorithm for this problem is simply
\begin{equation} \label{eq:newton1}
 a_j \leftarrow a_j + \frac{d_j e^{-a_j} - 2 \eta (a_j -c) - 1}{d_je^{-a_j} + 2\eta} 
\end{equation}
It can be readily shown that the solution to these $q$ optimization problems
produces  $\widehat{a}_{\eta,1} \ge \ldots \ge \widehat{a}_{\eta,q}$, with the inequalities being strict whenever the corresponding
inequalities for the corresponding sample eigenvalues are strict. 

For the shape version of this penalty, i.e.\ for $\Pi_{R,s}$, we have 
$\pi_{R,s}(a) = \sum_{j=1}^q (a_j - \overline{a})^2$ and hence $L_q(a;d;\eta) = \sum_{j=1}^q \{d_je^{-a_j} + a_j + \eta (a_j - \overline{a})^2\}$. For this case, the Newton-Raphson algorithm is given by
\begin{equation} \label{eq:newton2}
a_j \leftarrow a_j + \frac{g_j + \beta \sum_{k=1}^q g_k/\delta_k}{\delta_j}, 
\end{equation}
where  $g_j = d_j e^{-a_j} - 1 - 2 \eta (a_j - \overline{a})$,  $\delta_j = d_j e^{-a_j} + 2 \eta$
and $\beta =  2\eta q^{-1}/\{1 - 2\eta q^{-1} \sum_{i=1}^q \delta_j^{-1}\}$.

\section{Simulation study} \label{Sec:Sim}
In this section, we conduct a simulation study to compare the performance of the following five covariance estimators: \\

\begin{tabular}{rl}
S: & the sample covariance matrix, \\[2pt]
LogF: & the penalized covariance matrix proposed by \cite{Deng-Tsui:2013} with penalty $\|A\|^2_F$, 
\\ & where $A = \log \Sigma$,  \\[2pt]
sLogF: & our proposed shape penalized covariance matrix based on $\|A - \{\tr(A)/q\}I_q\|^2_F$, \\[2pt]
mLogF: & the adjusted penalized covariance matrix proposed by \cite{Yu:2017} based on \\ & $\|A - \widehat{m}I_q\|^2_F$, 
with $\widehat{m}$ being an estimate of $m(\Sigma) = tr(A)/q$, and \\[2pt]
dLogF: & an adjusted penalized covariance matrix based on $\|A - \log(\bar{d})I_q\|^2_F$, where $\bar{d} = \tr(S_n)/q$, 
\\ & i.e.\ the average of the sample eigenvalues. 
\end{tabular} 
\mbox{ } \\[8pt]
Comparisons of \emph{LogF} and \emph{mLogF} to other penalized covariance estimators are given in \cite{Deng-Tsui:2013}
and \cite{Yu:2017}. 

As the tuning constant $\eta \to \infty$, the estimator \emph{LogF} goes to the identity matrix and so one would anticipate 
its performance would be poor whenever $\bar{\lambda} = \tr{\Sigma}/q$ is far from one. This would be particularly problematic
when heavy tuning is needed, as would be the case whenever the roots of $\Sigma$ are not well separated or in general when 
$S_n$ is singular. As noted by \cite{Yu:2017}, this weakness can be alleviated by using the estimator \emph{mLogF}. Alternatively, the
estimators \emph{sLogF} or \emph{dLogF} can be considered. As shown in the appendix, the estimator \emph{sLogF} goes to 
$\bar{d}I_q$ as $\eta \to \infty$.
On the other hand, an adjusted estimator, i.e.\ one using a penalty of the form $\|A - c I_q\|^2_F$, 
goes to $e^c I_q$ as $\eta \to \infty$. Consequently, the estimators \emph{mLogF} and \emph{dLogF} go to 
$e^{\widehat{m}} I_q$  and $\bar{d}I_q$ respectively as $\eta \to \infty$.  

The performance of the estimator \emph{mLogF} depends on the definition of $\widehat{m}$. \cite{Yu:2017} observe that the 
simple choice $\widehat{m}_0 = m(S_n) = \tr(S_n)/q$ is known to underestimate $m(A)$. They propose using a bias corrected estimator of 
the form $\widehat{m}_1 = m(S_n) + b_{n,q}$ when $q < n$, and a Bayesian estimator for $\widehat{m}_3$ when $q \ge n$; 
see \cite{Yu:2017} for details. We use their proposed choices of $\widehat{m}$ in our simulation study. When using 
$\widehat{m}_1$ the estimator \emph{mLogF} shrinks the eigenvalues of $S_n$ towards 
$e^{\widehat{m}_1} \propto (\det S_n)^{1/q} = \bar{d}_g$, the geometric
mean of the eigenvalues of $S_n$. We surmise it would be better to shrink them towards the arithmetic mean 
since $\bar{d}$ is the minimum variance unbiased estimator of $\lambda$ when random sampling from a spherical multivariate
normal distribution with $\Sigma = \lambda I$.  In particular, we anticipate our proposed estimators \emph{sLogF} and 
\emph{dLogF}, which both shrink the eigenvalues of $S_n$ towards $\bar{d}$, will have a better performance in settings 
where heavy tuning is needed. The results of our simulation study, reported in Table \ref{Tab:Sim},
supports this heuristic argument.

For the simulations, we consider $q = 60$ dimensional data arising as a random sample from a multivariate normal 
distribution with mean $\mu = 0$ and covariance matrix $\Sigma = \{\sigma_{ij}\}$. The four different covariance 
models used in the simulations are listed below, along with the corresponding mean $\bar{\lambda}$ and 
standard deviation $s_{\lambda}$ of their eigenvalues. \\[-4pt]

\begin{tabular}{rl}
Model 1: & An MA(2) model where $\sigma_{ii}=10, \sigma_{i,i-1}=\sigma_{i-1,i}=0.1, \sigma_{i,i-2}=\sigma_{i-2,i}=0.05$,  
and \\ &$\sigma_{ij} = 0$ otherwise. Here $\bar{\lambda} = 10.0$ and $s_{\lambda} = 1.16$. \\[2pt]
Model 2: & An AR(1) model where $\sigma_{ij} = 0.5\rho^{|i-j|}$ and $\rho=0.3$. Here $\bar{\lambda} = 0.5$ 
and $s_{\lambda} = 0.22$. \\[2pt]
Model 3: & $\Sigma^{-1}=\{\sigma^{ij}\}$ where $\sigma^{ii}=1$ and $\sigma^{ij}=0.6$ for 
$i\neq j$. Here $\bar{\lambda} = 2.46$ and $s_{\lambda} = 0.32$.  \\[2pt]
Model 4: & $\Sigma=5I_q$. Here $\bar{\lambda} = 5.0$, and $s_{\lambda} = 0$. \\ & \\[-4pt]
\end{tabular}

To evaluate the performance of the different estimators under the various covariance models, four measures
of the discrepancy between the estimated covariance matrix and the true $\Sigma$ are computed. \\[-4pt]

\begin{tabular}{rl}
Fnorm: & $||\widehat{\Sigma}-\Sigma||_F=\sqrt{\sum_{i,j}(\widehat{\sigma}_{ij}-\sigma_{ij})^2})$. \\[2pt]
$L_1$: & $||\widehat{\Sigma}-\Sigma||_1=\max_j\sum_i|\widehat{\sigma}_{ij}-\sigma_{ij}|$. \\[2pt]
op-norm: & $||\widehat{\Sigma}-\Sigma||_{op}=\max_j|\tilde{\sigma}_j|$, where $\tilde{\sigma}_j$'s 
are the singular values of $\widehat{\Sigma}-\Sigma$. \\[2pt]
$\Delta_1$: & $|\widehat{\lambda}_1-\lambda_1|$, the absolute difference between the largest eigenvalues of 
$\widehat{\Sigma}$ and $\Sigma$, \\ & \\[-4pt]
\end{tabular}

We follow the simulation protocol used by both \cite{Deng-Tsui:2013} and  \cite{Yu:2017}.
For each covariance model, $2n$ data points are generated, with the first $n$ observations serving as a training set 
and last $n$ observations serving as a validation set. The tuning parameter for any particular method is selected to be
the value of $\eta$ which minimizes the non-penalized loss $l(\widehat{\Sigma}_{\eta,o}; S_{n,1})$ defined by \eqref{eq:lik},
where $\widehat{\Sigma}_{\eta,o}$ is the penalized covariance estimate based on the training set, and $S_{n,1}$ is
the sample covariance matrix of the validation set. Since the true mean $\mu = 0$ and interest lies in the performance of the
estimators of $\Sigma$, the non-centered sample covariance matrices $S_n = n^{-1}\sum_{i=1}^n x_i x_i^\tran$ are used in the 
simulations. We consider three values for the sample size $n=  30, 60, 120$, with the dimension being $q=60$ in each case. 
The simulations are repeated $100$ times and the means and standard deviations (in parenthesis) over the $100$ trials
for each of the discrepancy measures are reported in Table \ref{Tab:Sim}. The means and standard deviations over the 
$100$ trials of the value of the selected tuning parameter $\eta$ are also reported.

\begin{table}[!htbp] 
\vspace*{-1.5cm}
\caption{Simulation results for the performance of the five estimators of covariance under Models 1-4 based on four 
performance measures. Averages and standard deviation are calculated from 100 runs. \label{Tab:Sim}}
\vspace*{0.2cm}
 \centering
  \scriptsize
\begin{tabular}{ c c | c c c c c | c c c c c }
\hline
&&&&&&&&&&& \\[-6pt]
  &   & \multicolumn{5}{c|}{Model 1} & \multicolumn{5}{c}{Model 2} \\
$n$ & Method & Fnorm & $L_1$ & op-norm &$\Delta_1$  &$\eta$& Fnorm & $L_1$ & op-norm &$\Delta_1$ &$\eta$  \\
\hline
&&&&&&&&&&& \\[-6pt]
{n = 120} & sLogF  & 1.81*   & 1.09*   & 0.46*   & 0.23*   & 69.83   & 1.54   & 0.84*   & 0.46   & 0.28   & 2.21   \\
 &       & (0.58) & (0.68) & (0.17) & (0.16) & (37.19) & (0.03) & (0.05) & (0.01) & (0.03) & (0.34) \\
  & dLogF & 1.87&  1.13&  0.47&  0.23& 66.83& 1.46* & 0.99& 0.42* & 0.15* & 1.34\\
  &      & (0.62)&  (0.69)  &(0.20)  &(0.18)& (37.57)& (0.02)& (0.06) &(0.01)& (0.03)& (0.16)\\
   & mLogF    & 1.98   & 1.16   & 0.50   & 0.24   & 66.38   & 1.49   & 0.95   & 0.45   & 0.21   & 1.56   \\
   &       & (0.71) & (0.71) & (0.22) & (0.20) & (37.72) & (0.03) & (0.05) & (0.01) & (0.03) & (0.19)\\
    & LogF & 41.58  & 36.42  & 9.05   & 7.42   & 0.13    & 2.29   & 2.07   & 0.70   & 0.45   & 0.36   \\
    &       & (0.43) & (1.58) & (0.12) & (0.66) & (0.29)  & (0.04) & (0.08) & (0.03) & (0.04) & (0.01) \\
 & S     & 55.36  & 57.46  & 17.90  & 17.60  & NA      & 73.63  & 12.30  & 10.01  & 8.61   & NA     \\
&       & (1.26) & (3.14) & (1.17) & (1.17) & NA      & (0.06) & (0.15) & (0.02) & (0.08) & NA     \\
 \hline
&&&&&&&&&&& \\[-6pt]
{n = 60} & sLogF  & 2.08*   & 1.15*   & 0.52*   & 0.28*   & 78.29   & 1.64   & 0.72*   & 0.46   & 0.34   & 5.46   \\
&       & (0.81) & (0.61) & (0.19) & (0.18) & (30.90) & (0.02) & (0.04) & (0.01) & (0.02) & (0.62) \\
& dLogF & 2.22&  1.36&  0.56&  0.28& 73.48&1.57* & 0.94& 0.43* & 0.19* & 2.58\\
&      & (1.00)&  (0.95)&  (0.27)&  (0.25)& (35.06)& (0.02)& (0.08)& (0.01)& (0.04)& (0.49)\\
 & mLogF    & 60.71  & 20.02  & 9.66   & 3.95   & 1.10    & 3.31   & 1.35   & 0.81   & 0.54   & 0.90   \\
 &       & (0.45) & (0.64) & (0.04) & (0.28) & 0.00    & (0.04) & (0.04) & (0.01) & (0.03) & (0.05)\\
 & LogF & 58.03  & 20.80  & 9.45   & 3.18   & 1.10    & 2.76   & 2.53   & 0.88   & 0.61   & 0.50   \\
 &       & (0.21) & (0.69) & (0.02) & (0.28) & (0.00)  & (0.06) & (0.13) & (0.05) & (0.06) & (0.02) \\
& S     & 78.19  & 85.01  & 27.86  & 27.56  & NA      & 73.70  & 13.57  & 10.13  & 8.11   & NA     \\
 &       & (2.45) & (6.20) & (2.09) & (2.09) & NA      & (0.08) & (0.28) & (0.02) & (0.14) & NA     \\
\hline
&&&&&&&&&&& \\[-6pt]
{n = 30} & sLogF  & 2.81*   & 1.58*   & 0.66*   & 0.34*   & 78.30   & 1.69   & 0.59*   & 0.45   & 0.39   & 27.42   \\
&       & (1.27) & (0.93) & (0.25) & (0.25) & (30.92) & (0.02) & (0.08) & (0.01) & (0.04) & (24.70) \\
& dLogF & 2.99&  1.97&  0.81&  0.48& 72.15& 1.64* & 0.89& 0.44*& 0.22*& 5.16\\
&      &(1.39)&  (1.57)&  (0.52)&  (0.53)& (35.62) &(0.02) &(0.15) &(0.02) &(0.07) &(2.17)\\
& mLogF    & 63.81  & 15.80  & 9.17   & 5.76   & 4.18    & 3.12   & 1.14   & 0.77   & 0.59   & 4.10    \\
&       & (0.34) & (0.48) & (0.03) & (0.24) & (0.27)  & (0.04) & (0.03) & (0.01) & (0.02) & (0.35) \\
& LogF & 65.40  & 15.75  & 9.35   & 6.07   & 4.13    & 3.20   & 2.40   & 0.94   & 0.61   & 1.10    \\
&       & (0.15) & (0.42) & (0.01) & (0.19) & (0.17)  & (0.05) & (0.11) & (0.04) & (0.06) & (0.00)  \\
& S     & 110.65 & 125.53 & 44.51  & 44.21  & NA      & 73.79  & 15.41  & 10.22  & 7.30   & NA      \\
&       & (3.79) & (9.51) & (3.58) & (3.58) & NA      & (0.11) & (0.44) & (0.01) & (0.24) & NA      \\
\hline 
\mbox{ } \\
%

\hline
&&&&&&&&&&& \\[-6pt]
  &   & \multicolumn{5}{c|}{Model 3} & \multicolumn{5}{c}{Model 4} \\
$n$ & Method & Fnorm & $L_1$ & op-norm & $\Delta_1$ & $\eta$& Fnorm & $L_1$ & op-norm &$\Delta_1$ &$\eta$  \\
\hline
&&&&&&&&&&& \\[-6pt]
{n=120} & sLogF  & 2.45*   & 2.43   & 2.31   & 0.08*  & 16.05  & 0.61*   & 0.37*   & 0.14*   & 0.12*   & 124.21  \\
 &       & (0.04) & (0.04) & (0.05) & (0.07) & (4.89) & (0.40) & (0.38) & (0.11) & (0.11) & (83.72) \\
 & dLogF & 2.45*&  2.40*&  2.33&  0.15& 13.48&0.65&   0.40&   0.15&   0.14& 117.74\\
 &      &(0.04)& (0.03)& (0.06)& (0.10)& (5.71)&(0.42)&  (0.39)&  (0.13)&  (0.13)& (84.25)\\
 & mLogF   & 2.61   & 2.51   & 2.21   & 0.08   & 13.91  & 0.72   & 0.41   & 0.17   & 0.16   & 118.41  \\
 &       & (0.10) & (0.05) & (0.05) & (0.06) & (5.42) & (0.46) & (0.40) & (0.14) & (0.15) & (85.33)\\
  & LogF & 9.03   & 5.12   & 1.82*   & 0.16   & 0.83   & 26.31  & 27.04  & 8.26   & 8.26   & 0.01    \\
  &       & (0.12) & (0.23) & (0.02) & (0.12) & (0.08) & (0.57) & (1.47) & (0.55) & (0.55) & (0.00)  \\
 & S     & 59.99  & 21.15  & 10.28  & 3.37   & NA     & 47.52  & 32.02  & 9.54   & 3.65   & NA      \\
 &       & (0.23) & (0.66) & (0.00) & (0.29) & NA     & (0.32) & (1.33) & (0.07) & (0.59) & NA      \\
 \hline
&&&&&&&&&&& \\[-6pt]
{n=60} & sLogF  & 2.67   & 2.51   & 2.20   & 0.16*   & 10.00  & 0.79*   & 0.45*   & 0.17*   & 0.12*   & 78.54   \\
&       & (0.13) & (0.05) & (0.05) & (0.06) & 0.00   & (0.51) & (0.32) & (0.11) & (0.12) & (30.98) \\
& dLogF & 2.54*& 2.43*& 2.31& 0.28& 9.80&0.87&  0.55&  0.21&  0.18& 73.77\\
&     &(0.06)& (0.04)& (0.05)& (0.07)& (0.78)&(0.60)  &(0.49)  &(0.15)  &(0.16)& (34.94)\\
 & mLogF   & 14.78  & 5.60   & 2.41   & 0.74   & 0.89   & 30.35  & 10.00  & 4.79   & 1.82   & 1.10    \\
 &       & (0.18) & (0.19) & (0.01) & (0.11) & (0.06) & (0.22) & (0.31) & (0.01) & (0.14) & (0.00) \\
 & LogF & 10.45  & 4.20   & 1.66*   & 0.54   & 2.81   & 26.94  & 8.99   & 4.24   & 1.58   & 1.92    \\
 &       & (0.07) & (0.12) & (0.01) & (0.07) & (0.24) & (0.61) & (0.98) & (0.06) & (0.46) & (0.39)  \\
 & S     & 61.57  & 27.13  & 10.29  & 0.93   & NA     & 55.05  & 44.49  & 10.16  & 8.63   & NA      \\
 &       & (0.27) & (1.07) & (0.00) & (0.49) & NA     & (0.55) & (2.28) & (0.28) & (1.04) & NA      \\
\hline
&&&&&&&&&&& \\[-6pt]
{n=30} & sLogF  & 2.53*   & 2.49*   & 2.37   & 0.09*   & 173.79   & 1.15*   & 0.57*   & 0.22*   & 0.18*   & 137.38  \\
 &       & (0.10) & (0.04) & (0.08) & (0.07) & (125.83) & (0.71) & (0.53) & (0.14) & (0.15) & (75.73) \\
 & dLogF & 2.56&   2.50&   2.39&   0.16& 134.08&1.30&   0.78&   0.32&   0.29& 129.55\\
 &       &(0.13)   &(0.09)  &(0.09)   &(0.18) &(127.69)&(0.81)  &(0.88)  &(0.30)  &(0.31)& (82.46)\\
 & mLogF    & 14.65  & 4.57   & 2.08   & 1.18   & 4.44     & 30.71  & 7.85   & 4.32   & 2.52   & 4.30    \\
 &       & (0.12) & (0.10) & (0.01) & (0.08) & (0.48)   & (0.21) & (0.29) & (0.02) & (0.15) & (0.40) \\
 & LogF & 11.04  & 3.93   & 1.57*   & 0.82   & 6.90     & 29.10  & 7.62   & 4.10   & 2.31   & 4.67    \\
 &       & (0.05) & (0.08) & (0.01) & (0.10) & (1.06)   & (0.14) & (0.37) & (0.01) & (0.18) & (0.52)  \\
& S     & 64.51  & 36.10  & 10.29  & 3.21   & NA       & 67.55  & 62.59  & 17.26  & 16.96  & NA      \\
&       & (0.27) & (1.89) & (0.00) & (0.90) & NA       & (1.07) & (3.77) & (1.79) & (1.79) & NA      \\
\hline
\end{tabular}
\end{table}

The five estimators in Table \ref{Tab:Sim} are listed by the order of the overall performance over the 
various models. The estimator which performed the best for a given sample size and given discrepancy measure
is noted with an asterisk (*). For every discrepancy measure, our proposed \emph{sLogF} outperforms all the other 
estimators under models 1 and 4 at every sample size. Under models 2 and 3, with one exception, either \emph{sLogF} or 
\emph{dLogF} is the best performing estimator, depending on the particular discrepancy measure used. The
notable exception is under model 3, for which \emph{LogF} performs best under the operator norm. Overall, the 
performance of our two proposed estimators \emph{sLogF} and \emph{dLogF} are similar for all four models. The 
performance of \emph{mLog} is also similar to these two estimators when $n = 120$, but performs considerably worse when 
$n = 60$ or $n=30$. As previously surmised, \emph{sLogF} and \emph{dLogF} have particularly better performance than the
other estimators whenever their tuning parameters tend to be large. Finally, as suspected, the sample covariance
matrix uniformly performs the worse.


\section{An example: Sonar data} \label{Sec:Ex}
As an example, we consider the sonar data set obtained from University of California Irvine Machine Learning Repository, which
was developed and first analyzed by \cite{Gorman:1988} This data set consists of $208$ multivariate observations of dimension
$q = 60$. For each observation, the $60$ variables correspond to the average energy over a particular frequency 
band obtained by bouncing sonar signals off of an object under various conditions, with $111$ observations labeled 
M (metal cylinder) and the other $97$ observations labeled R (rock). 

Our goal here is to study the relative performance of covariance estimators when used within linear discriminant 
analysis (LDA) to classify an observation as either M or R.  As in the simulation study, the data set is randomly partitioned 
into a training set of size $78$ for estimating the covariance matrix, a validation set of size $78$ for selecting the 
tuning parameter and a test set of size $52$ for computing the misclassification error. The covariance estimators
being compared are those considered in the simulation study in section \ref{Sec:Sim}. Here, though, the estimators 
are based on the pooled sample covariance matrix of the two groups M and R.

\begin{figure}[!hbtp]
		\centering
		\includegraphics[width=12cm]{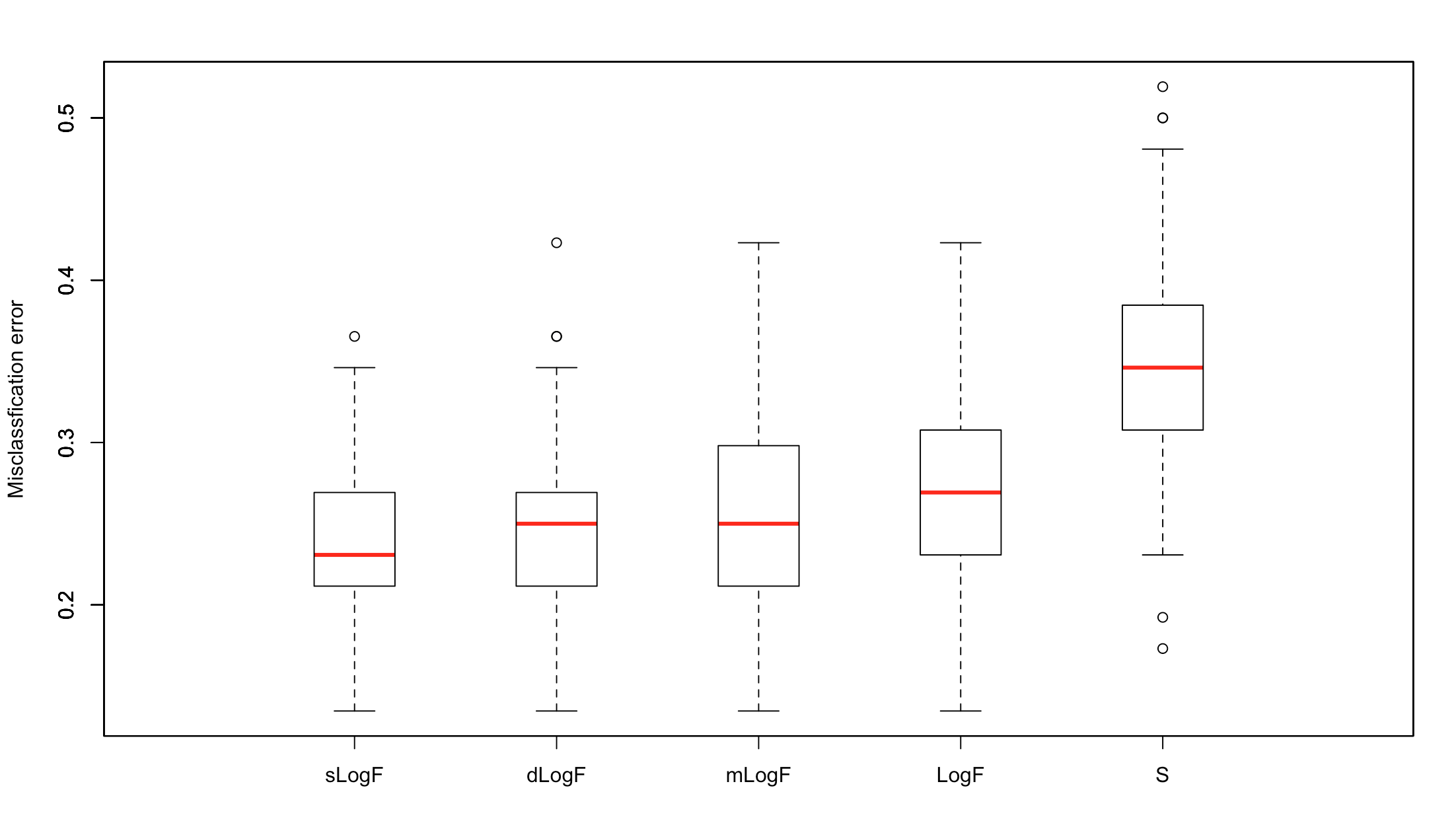}
		\caption{Boxplot of misclassification errors over 100 runs.}
		\label{fig:sonar}
	\end{figure}

The above procedure is independently repeated for 100 times.  A boxplot of the misclassification errors over these $100$
trials are displayed in Figure \ref{fig:sonar}, and the mean and standard deviation of the misclassification errors 
are showed in Table \ref{tab:sonar}. Finally, Table \ref{tab:sonar-freq} displays the frequency over the $100$ trials 
that a given estimator (row) has a lower classification rate than another estimator (column). For example,
\emph{sLogF} has less misclassification errors than \emph{mLogF} in $46$ of the $100$ runs, and more misclassification errors
in $22$ of the runs, with the two estimators having the same misclassification rate in the other $32$ runs.
Among the estimators of the covariance matrix considered here, our proposed \emph{sLogF} estimator performs best.
	
	\begin{table}[!htbp]
		\centering
		\caption{Means and standard deviations of the misclassification error over 100 runs.}
		\vspace{0.2cm}
		\centering
		\begin{tabular}{ |c  | c |c |c|c|c|}
			\hline
			  &sLogF& dLogF &mLogF&LogF&S\\
			\hline
			 $\qquad$Mean$\qquad$ &0.239 &0.247&0.259 &0.267  &0.349\\
			\hline
			 S.D. &0.048 &0.051&0.062&0.059   &0.071\\
			\hline
		\end{tabular}
		\label{tab:sonar}
	\end{table}
	
	\begin{table}[!htbp]
		\centering
		\caption{Frequency of less misclassifications using row estimator versus column estimator out of 100 runs.}
		\vspace{0.2cm}
		
		\begin{tabular}{|c|c|c|c|c|c|}
			\hline
			& sLogF & dLogF & mLogF & LogF & $\quad$ S$\quad $ \\
			\hline
			$\quad$sLogF$\quad$ & 0     & 38    & 46    & 60   & 89 \\
			\hline
			dLogF & 28    & 0     & 31    & 56   & 93 \\
			\hline
			mLogF & 22    & 18    & 0     & 54   & 84 \\
			\hline
			LogF  & 26    & 31    & 35    & 0    & 83 \\
			\hline
			S     & 5     & 6     & 12    & 14   & 0 \\
			\hline
		\end{tabular}
		\label{tab:sonar-freq}
	\end{table}
\newpage
\section{Appendix: Proofs and some technical details}
\subsection*{\small \emph{Counterexamples to the equivalency of g-convexity and convexity in log.}}
Lemma 1.14 in \cite{Wiesel-Zhang:2015} states that $x^\tran \Sigma^{-1} x$ is a strictly g-convex 
function of $\Sigma$, which implies that $\tr(\Sigma^{-1}S_n)$ is g-convex for $S_n \ne 0$. 
It is difficult to show analytically whether or not $\tr(\Sigma^{-1}S_n)$ is convex in $\log \Sigma$ for 
a given $S_n$, and almost all randomly generated counterexamples tend to imply that it true. After extensive
trials, though, the following counterexample was found which shows that $\{\Sigma^{-1}\}_{11}$ is not a convex 
function of $\log \Sigma$, and consequently $\tr(\Sigma^{-1}S_n)$ cannot be convex in $\log \Sigma$ in general.
For $q = 2$, let $A = \log \Sigma$ and choose 
\[ A_0 = \left[ \begin{array}{cc} 0 & -1 \\ -1 & 300 \end{array} \right] \quad \mbox{and} \quad  
A_1  = - \left[ \begin{array}{cc} 0 & 0.01 \\ 0.01 & 0.01 \end{array} \right].\] 
This gives
$\{e^{-(0.5 A_1 +0.5 A_2)}\}_{11} = 1.001690296 > 1.001688939 = 0.5 \{e^{-A_0}\}_{11} + 0.5 \{e^{-A_1}\}_{11}$,
and so $\{e^{-A}\}_{11} = \{\Sigma^{-1}\}_{11}$ is not convex in $A$. 

On the other hand, a function may be convex in $\log \Sigma$ but not g-convex in $\Sigma$. 
For example, the matrix $L_1$ norm on the elements of $\log \Sigma$, i.e.\  
$H(\Sigma) = \max_{1 \le k \le q} \sum_{j=1}^q | \{ \log \Sigma \}_{jk} |$, is convex in $\log \Sigma$.
The following counter-example, though, shows that it is not g-convex. For $q = 3$, choose
\[ \Sigma_0 = \left[ \begin{array}{ccc} 
1.00 & 0.30 & 0.09 \\ 0.30 & 1.00  & 0.30 \\  0.09 & 0.30 & 1.00 \end{array} \right] \quad \mbox{and} \quad  
\Sigma_1  = \left[ \begin{array}{ccc} 
 1.00 & 0.90 & 0.81 \\ 0.90 & 1.00  & 0.90 \\  0.81 & 0.90 & 1.00 \end{array} \right].\] 
This gives
$H(\Sigma_{0.5}) = 2.289438 > 2.284073 = 0.5~H(\Sigma_0) + 0.5~H(\Sigma_1)$,
and so $H(\Sigma)$ is not g-convex. 
 
\subsection*{\small \emph{Proof of Lemma \ref{Lem:sing}}}
The lemma follows if $L(\Sigma;S_n,\eta)$ is g-coercive. 
Consider any sequence in $\Sigma$ such that $\| \log \Sigma \|_F  \to \infty$.  Divide the proof into the following three cases:
a) $\log \det \Sigma \to \infty$, b) $|\log \det \Sigma|$ is bounded above, and c) $\log \det \Sigma \to -\infty$. 
For case (a), the result holds since both $\tr(\Sigma^{-1}S_n) \ge 0$ and $\Pi(\Sigma) \ge 0$.  For case (b), the result follows from
condition (i) since $\tr(\Sigma^{-1}S_n) \ge 0$ and $\Pi(\Sigma) \to \infty$. 

When case (c) holds, consider the two sub-cases: c1) $\lambda_1$ is bounded away from zero, and c2) $\lambda_1 \to 0$. If (c1) holds, condition (ii) 
implies $(\log \det \Sigma)/\Pi(\Sigma) \to 0$ and so  $\Pi(\Sigma) \to \infty$. Hence, for $\eta > 0$,
\[L(\Sigma;S_n,\eta)  = \tr(\Sigma^{-1}S_n) + \Pi(\Sigma)\{ (\log \det \Sigma)/\Pi(\Sigma) + \eta\} \to \infty. \]
If (c2) holds, since $\tr(\Sigma^{-1}S_n) \ge \tr(S_n)/\lambda_1$ and $\log \det \Sigma \ge \log \lambda_1 - (q-1)\log \lambda_q$, it follows that
\[L(\Sigma:S_n,\eta) \ge \tr(S_n)/\lambda_1 + q\log \lambda_1 + (q-1)\log(\lambda_q/\lambda_1) + \eta \Pi(\Sigma),\]
with $\tr(S_n)/\lambda_1 + q\log \lambda_1 \to \infty$. So, if $\lambda_1/\lambda_q \to 1$, then $L(\Sigma:S_n,\eta) \to \infty$.
Whereas, if $\lambda_1/\lambda_q$ is bounded away from one, then by condition (iii), $(q-1)\log(\lambda_q/\lambda_1) + \eta \Pi(\Sigma) =
\Pi(\Sigma)\{(q-1)\log(\lambda_q/\lambda_1)/\Pi(\Sigma) + \eta)$ is bounded below and so $L(\Sigma:S_n,\eta) \to \infty$.

\subsection*{\small \emph{Proof of Theorem \ref{Thm:g-log}}}

First, we show $(i)\Rightarrow (iii)$. Suppose that $F(\Sigma)$ is (strictly) g-convex, then by Lemma 3.6 of \cite{Duembgen-Tyler:2016}, 
$F(BD(e^x)B^\tran)$ is (strictly) convex in $x\in\R^q\backslash\{0\}$ for any non-singular $B$ of order $q$. Here, for $y \in \R^q$, $D(y)$ represents the
diagonal matrix with the elements of $y$ corresponding to its diagonal elements. Thus, by Lemma \ref{Lem:orth}, $f(x) = F(D(e^x))$ is (strictly) convex. 

Next, we show $(iii)\Rightarrow (i)$. Here, the concept of majorization plays an important role.  For a vector $v \in \R^q$, denote its ordered
values by $v_{(1)} \ge \cdots \ge v_{(q)}$. A vector $y \in \R^q$ is then said to majorize a vector $x \in \R^q$, denoted $x \prec y$ if and only if 
$\sum_{j=1}^k x_{(j)} \le \sum_{j=1}^k y_{(j)}$, with equality when $k = q$. As stated in Theorem 1.3 of \cite{Ando:1989}, $x \prec y$ if an only
if $x$ is a convex combination of coordinate permutations of $y$, i.e.\
\begin{equation} \label{eq:Ando}
x \prec y \Leftrightarrow x = \sum_{j=1}^q w_j P_j y, 
\end{equation} 
where, for $j = 1, \ldots, q$, $P_j$ is a permutation matrix of order $q$, hence orthogonal, and $w_j \ge 0$ with $\sum_{j=1}^q w_j = 1$. 
As a side note, the Birkhoff-von Neumann Theorem notes that $Q$ is a doubly stochastic matrix of order $q$ if and only if it has the
representation $Q = \sum_{j=1}^q w_j P_j$.
For $\Sigma > 0$, let $\lambda(\Sigma) = (\lambda_1(\Sigma), \ldots, \lambda_q(\Sigma))$ denote the vector of the ordered
eigenvalues of $\Sigma$. An important result given by Lemma 2.17 in \cite{Sra-Hosseini:2015} states 
\begin{equation} \label{log-maj}
\log(\lambda(\Sigma_t))\prec(1-t)\log(\lambda(\Sigma_0))+t\log(\lambda(\Sigma_1)),
\end{equation}
where $\Sigma_t$ is the geodesic curve from $\Sigma_0$ and $\Sigma_1$. So, by \eqref{eq:Ando}, we can express
\begin{equation} \label{eq:Q}
 \log\lambda(\Sigma_t)=Q\{(1-t)\log\lambda(\Sigma_0)+t\log\lambda(\Sigma_1)\},
\end{equation}
with $Q =\sum_{j=1}^q w_jP_j$ being defined as in \eqref{eq:Ando}.  Thus, 
\begin{align*}
    F(\Sigma_t) &=f(\log\lambda(\Sigma_t))=f\left(Q[(1-t)\log\lambda(\Sigma_0)+t\log\lambda(\Sigma_1)]\right) \\[2pt]
    &\le(1-t)f(Q\log\lambda(\Sigma_0))+tf(Q\log\lambda(\Sigma_1))  \\
    &\le(1-t)\sum_{j=1}^q w_j f(P_j\log\lambda(\Sigma_0))+t\sum_{j=1}^q w_j f(P_j\log\lambda(\Sigma_1)) \\
    &=(1-t)\sum_{j=1}^q w_j f(\log\lambda(\Sigma_0))+t\sum_{j=1}^q w_j f(\log\lambda(\Sigma_1)) \\
    &=(1-t)F(\Sigma_0)+tF(\Sigma_1). 
\end{align*}
The two inequalities above follow from condition (iii), i.e.\ $f$ is convex. Suppose now that $f$ is strictly convex, then the 
first inequality is strict unless $\lambda(\Sigma_0) = \lambda(\Sigma_1)$, and the second inequality is strict unless $Q = I$.  
Thus, both equality holds if and only if $\lambda(\Sigma_t)= \lambda(\Sigma_0)$ for $0 \le t \le 1$. However, 
since $\tr \Sigma$ is strictly g-convex, see e.g.\ Lemma 1.15 in \cite{Wiesel-Zhang:2015}, it follows that 
$\tr(\Sigma_t) < (1-t) \tr(\Sigma_0) + t ~ \tr(\Sigma_1) = \tr(\Sigma_0)$ for $0 < t < 1$, unless
$\Sigma_0 = \Sigma_1$. Thus, $F(\Sigma)$ is strictly g-convex.

Finally, we note the statement $(ii)\Leftrightarrow (iii)$  follows from the main theorem in \cite{Davis:1957}, 
at least in the convex case. The strictly convex case can be shown to hold by applying arguments analogous to 
those used in the $(i)\Leftrightarrow (iii)$ case.

\subsection*{{\small \emph{Proof that} {\boldmath $\Pi_{R,s}$} \emph{satisfies the conditions of Lemma \ref{Lem:sing}}.}}
Again let $a_j = \log \lambda_j$, and so $\overline{a} = q^{-1} \log \det \Sigma$, $\sum_{j=1}^q a_j^2 = \| \log \Sigma \|_F^2$
and $ss(a) \equiv \sum_{j=1}^q (a_j - \overline{a})^2  = \Pi(\Sigma)$. Condition (i) states 
that if $|\overline{a}|$ is bounded above and $\sum_{j=1}^q a_j^2 \to \infty$, then $ss(a) \to \infty$, which 
holds since $ss(a) =  \sum_{j=1}^q a_j^2 - q \overline{a}^2$.  Condition (ii) states that if $\overline{a} \to -\infty$ 
and $a_1$ is bounded below, then $\overline{a}/ss_a \to 0$. To show this, express 
$\overline{a}/ss(a) = \{ \overline{a} ss(b) \}^{-1}$, where $b_j = a_j/\overline{a}$. Since $a_1$ is bounded below,
$b_1 \to 0$ and $\sum_{j=1}^q b_j \to q$. Hence, $ss(b)$ must be bounded away from zero, which implies
$\overline{a}/ss(a) \to 0$. Condition (iii) states that if $a_1 \to -\infty$ and $a_1 - a_q > \epsilon > 0$,
then $(a_1 - a_q)/ss(a)$ is bounded above. This follows since $ss(a) \ge (a_1 - a_q)^2$ and so 
$(a_1 - a_q)/ss(a) \le 1/(a_1 - a_q) \le 1/\epsilon$.

\subsection*{\small \emph{Proof of Theorem \ref{Thm:shape}}} 
i) Let $\Sigma_0 \#_t \Sigma_1 := \Sigma_0^{1/2}(\Sigma_o^{-1/2}\Sigma_1 \Sigma_o^{-1/2})^t\Sigma_o^{1/2}, t\in[0,1]$,
and so $\Sigma_t =  \Sigma_0 \#_t \Sigma_1$ \citep{Sra-Hosseini:2015}. It readily follows that $V(\Sigma_t)=V(\Sigma_0)\#_tV(\Sigma_1)$, and so
\begin{align*}
    \Pi_s(\Sigma_t) &=\Pi\{V(\Sigma_t)\} =\Pi\{V(\Sigma_0)\#_tV(\Sigma_1)\}\\
    &\le (1-t)~\Pi\{V(\Sigma_0)\}+t~\Pi\{V(\Sigma_1)\} =(1-t)~\Pi_s(\Sigma_0)+t~\Pi_s(\Sigma_1). 
\end{align*}
ii) Let $A=\log\Sigma$, and define $G(A) = \Pi(e^A)$ and $G_s(A) = \Pi_s(e^A)$. The goal is to show that if $G(A)$ is convex in $A$, then 
$G_s(A)$ is also convex in $A$. Since $G_s(A) = G(\tilde{A})$, where $\tilde{A} \equiv A - (\tr(A)/q)*I$, and so  
\begin{align*}
    G_s&((1-t)A_0+tA_1) = G((1-t)\tilde{A}_0  + t\tilde{A}_1)  \\ 
    &\le (1-t)G(\tilde{A}_0) + t G(\tilde{A}_1) = (1-t)G_s(A_0)+tG_s(A_1.) 
\end{align*}

\subsection*{\small Proof of Theorem \ref{Thm:diag}.} 
The proof relies on the following well known extremal property of eigenvalues of symmetric matrices. 
Let $B$ be a symmetric matrix of order $q$, and let $C = [c_1 \cdots c_m]$ be of order $q \times m$, $m \le q$,
with orthonormal columns. Then $\tr\{C^\tran B C \} = \sum_{j=1}^m c_j^\tran B c_j$ is bounded above and below 
by the sum of the largest $m$  and the sum of the smallest $m$ eigenvalues of $B$ respectively.

Expressing $\Sigma = P\Lambda P^\tran$ in terms of its spectral value decomposition, let
$H = [h_1 \cdots h_q] = P_n^\tran P$, which is itself an orthogonal matrix. Define $\kappa_1 = 1/\lambda_1$ and 
$\kappa_j = 1/\lambda_j - 1/\lambda_{j-1}$ for $j \ne 1$. Inverting this relationship gives 
$\lambda_j^{-1} = \sum_{k=1}^j \kappa_k$. Since $\kappa_j  \ge 0$, the above noted extremal 
property of eigenvalues of a symmetric matrix implies
\[ \begin{array}{lcl}
 \tr\{\Sigma^{-1}S_n\} &=& \tr\{\Lambda^{-1}H^\tran D_n H\}  = \sum_{j=1}^q \lambda_j^{-1} h_j^\tran D_n h_j = 
\sum_{k=1}^q \kappa_k \left\{\sum_{j=k}^q h_j^\tran D_n h_j \right\} \\
 &\ge&  \sum_{j=1}^q \kappa_j \left\{\sum_{k=j}^q d_k \right\} = \sum_{j=1}^q d_j/\lambda_j = \tr\{\Lambda^{-1}D_n\}, 
\end{array} \]
with equality when $Q = P$.  
The lemma follows since $\det \Sigma = \det \Lambda$ and $\Pi(\Sigma) = \Pi(\Lambda)$. 

\subsection*{{\small \emph{Limiting behavior of the sLogF estimator as} {\boldmath $\eta \to \infty. $}}}
As $\eta \to \infty$, the penalty term $\Pi_{R,s}(\widehat{\Sigma}_\eta)$ must go to $0$, which 
implies $\widehat{\Sigma}_\eta$ is proportional to $I_q$ in the limit.
The eigenvalues of $\widehat{\Sigma}_\eta$ correspond to the unique critical point of  
$\sum_{j=1}^q \{d_j e^{-a_j} + a_j + \eta (a_j -\overline{a})^2\}$, where again $a_j = \log \lambda_j$,
which in turn corresponds to the unique solution to the set of equations 
$d_j e^{-a_j} = 1 + 2 \eta (a_j -\overline{a})$ for $j = 1, \ldots, q$. By taking the sum, we obtain 
$q =  \sum_{j=1}^q d_j e^{-a_j} = \sum_{j=1}^q d_j/\lambda_j$ for any $\eta \ge 0$. Hence, since the eigenvalues of 
$\widehat{\Sigma}_\eta$ approach each other as $\eta \to \infty$, it follows that $\widehat{\lambda}_j \to \overline{d}$
or  $\widehat{\Sigma}_\eta \to \overline{d} I_q$.

\end{document}